\documentclass[12pt]{amsart}
\usepackage{amsmath,amssymb}

\newcommand{\Ch}{\mathrm{Ch}}

\newcommand{\mcl}{\mathrm{cl}}

\newcommand{\rr}{\mathrm{r}}

\newtheorem{theorem}{Theorem}[section]

\theoremstyle{definition}

\newtheorem{example}[theorem]{Example}

\theoremstyle{remark}

\numberwithin{equation}{section}

\begin{document}
\title[]
{Polynomially spectrum-preserving maps between commutative 
Banach 
algebras}

\author[O.~Hatori]{Osamu~Hatori}
\address{Department of Mathematics, Faculty of Science, 
Niigata University, Niigata 950-2181 Japan}
\curraddr{}
\email{hatori@math.sc.niigata-u.ac.jp}
\thanks{The authors were partly 
supported by the Grants-in-Aid for Scientific 
Research, The 
Ministry of Education, Science, Sports and Culture, Japan.}

\author[T.~Miura]{Takeshi Miura}
\address{Department of Basic Technology,
Applied Mathematics and Phisics, 
Yamagata University, Yonezawa 992-8510 Japan}
\curraddr{}
\email{miura@yz.yamagata-u.ac.jp}

\author[H.~Takagi]{Hiroyuki Takagi}
\address{Department of Mathematical Sciences,
Faculty of Science, Shinshu University, Matsumoto 390-8621 Japan}
\curraddr{}
\email{takagi@math.shinshu-u.ac.jp}

\date{}

\begin{abstract}
Let $A$ and $B$ be unital semi-simple commutative Banach 
algebras. In this paper we study two-variable polynomials $p$ which satisfy 
the following property: a map $T$ from $A$ onto $B$ such that the 
equality
\[
\sigma (p(Tf,Tg))=\sigma (p(f,g)), \quad f,g \in A
\]
holds is an algebra isomorphism.
\end{abstract}

\maketitle

\section{Introduction}
The study of spectrum-preserving linear maps between Banach algebras 
dates back to Frobenius \cite{fro} who studied linear maps on matrix 
algebras which preserve the determinant. After over 100 years spectrum-
preserving maps  are studied for Banach algebras and
the following conjecture seems to be 
still open: 
any 
spectrum-preserving linear map 
from a unital Banach algebra onto a unital semi-simple  
Banach algebra that preserves the unit is a Jordan morphism. 
The Gleason, Kahane and \.Zelazko theorem \cite{gl,kz,ze} 
asserts that a unital linear 
functional defined on a Banach algebra is multiplicative if 
it is invertibility preserving and the theorem 
has inspired a number of papers on 
the subjects. For commutative Banach algebras 
it is a straightforward 
conclusion of the theorem of Gleason, Kahane and \.Zelazko 
that a unital and 
spectrum-preserving linear map from a Banach algebra into 
a semi-simple {\it commutative} Banach algebra is a homomorphism. 
Thus the problems on spectrum-preserving linear maps 
 mainly concerns with 
non-commutative Banach algebras and 
has seen much progress 
recently \cite{au1,js,ne,so}.

Without assuming linearity, non-multiplicative and invertibility 
preserving maps are almost arbitrary, and spectrum-preserving 
maps which are not linear nor multiplicative are also 
possible even in the case of commutative Banach algebras. 
On the other hand, spectrum-preserving maps on Banach algebras 
which are not assumed to be linear are studied by several authors 
\cite{hmt1,hmt2,h,l,lt,ps,ra,rr1,rr2} recently. 
In this paper we study linearity and 
multiplicativity of spectrum-preserving 
maps between commutative Banach algebras 
under additional assumptions.

Let $A$ and $B$ be unital Banach algebras. Suppose that $S$ is an algebra 
isomorphism from $A$ onto $B$. Then we have that the equality 
\[
\sigma (p(Tf))=\sigma (p(f)), \quad f \in A
\]
holds for every polynomial $p$, where $\sigma (\cdot )$ denotes the 
spectrum. But the converse does not hold in general. Suppose that 
$X$ is a compact Hausdorff space and $C(X)$ denotes the algebra of all 
complex-valued continuous functions on $X$. For each $f\in C(X)$, 
$\pi_{f}$ denotes a self homeomorphism on $X$. Put a map $T$ from 
$C(X)$ into itself by 
\[
Tf=f\circ\pi_f
\]
for every $f\in C(X)$. 
Then $T$ may not be linear nor multiplicative while
\[
\sigma (p(Tf))=\sigma (p(f)), \quad f\in C(X)
\]
holds for every polynomial. But the situation is very different for 
polynomials of two variables. In this paper we show that
for certain two-variable polynomials 
$p(z,w)$ the following holds: 
a map $T$ from a unital semi-simple commutative Banach algebra $A$
onto another one $B$ 
is an algebra isomorphism if the equation
\[
\sigma (p(Tf,Tg))=\sigma (p(f,g)), \quad f,g \in A
\]
holds. 

\section{preliminary}
Let $X$ be a compact Hausdorff space. The algebra of all complex-valued 
continuous functions on $X$ is denoted by $C(X)$. For a subset $K$ of $X$ 
the uniform norm on $K$ is denoted by 
$\|\cdot \|_{\infty (K)}$. A uniform algebra 
on $X$ is a uniformly closed subalgebra of $C(X)$ which separates the 
points of $X$ and contains the constant functions. For a uniform algebra 
$A$ on $X$, $P(A)$ denotes the set of all peaking functions in $A$. 
The set of all weak peak points for $A$ is the Choquet boundary and 
denoted by $\Ch (A)$. See \cite{B,ga} for theory of uniform algebras. 
Let ${\mathcal A}$ be a commutative Banach algebra.
We denote the maximal ideal space of ${\mathcal A}$ by $M_{\mathcal A}$ and 
the Gelfand transformation of $f\in {\mathcal A}$ is denoted by 
$\hat f$. The spectral radius for $f\in {\mathcal A}$ 
is denoted by $\rr (f)$ and the spectrum of $f$ is denoted by $\sigma (f)$.
The complex number field is denoted by ${\mathbb C}$.

\section{A conclusion of a theorem of Kowalski and S\l odkowski}
Kowalski and S\l odkowski \cite{ks} proved the following surprising generalization of a theorem of Gleason, Kahane and \.Zelazko.
\begin{theorem}
Let $A$ be a Banach algebra and $\phi$ a complex-valued map defined on $A$. 
Suppose that 
\[
\varphi (f) -\varphi (g) \in \sigma (f-g)
\]
holds for every pair $f$ and $g$ in $A$. Then $\varphi-\varphi (0)$ is 
linear and multiplicative.
\end{theorem}

Applying the above theorem we see the following. 

\begin{theorem}\label{cks}
Let $A$ be a Banach algebra and $B$ a semi-simple commutative Banach algebra, 
and $p(z,w)=az+bw$ ($ab\ne 0$). Suppose that $T$ is a 
(not necessarily linear) map from $A$ into $B$ which satisfies that the 
inclusion 
\[
\sigma (p(Tf,Tg)\subset \sigma (p(f,g))
\]
holds for every pair $f$ and $g$ in $A$. Then we have the following.

(1) If $a+b\ne 0$, then $T$ is linear and multiplicative. 

(2) If $a+b=0$, then $T-T(0)$ is linear and multiplicative. 
\end{theorem}

\begin{proof}
First we show that 
\[
\sigma (Tf-Tg)\subset \sigma (f-g), \quad f,g\in A
\]
holds. Let $f,g\in A$. Since $a\ne 0$, we have
\[
\sigma (Tf+\frac{b}{a}Tg)\subset \sigma (f+ \frac{b}{a}g), 
\]
so
\[
\sigma(T(-\frac{b}{a}g)+\frac{b}{a}Tg)\subset \sigma (-\frac{b}{a}g+
\frac{b}{a}g)=\{0\}
\]
by putting $f=-\frac{b}{a}g$. Thus the equality
\[
T(-\frac{b}{a}g)=-\frac{b}{a}Tg
\]
holds for every $g\in A$ since $B$ is semi-simple. It follows that 
\begin{multline*}
\sigma(Tf-Tg)=\sigma (Tf-T(-\frac{b}{a}(-\frac{a}{b}g))
\\=
\sigma (Tf+\frac{b}{a}T(-\frac{a}{b}g))\subset \sigma (f-g)
\end{multline*}
holds for every pair $f$ and $g$ in $A$. 

Put a map $S$ from $A$ into $B$ by $Sf=Tf-T(0)$. Then $S$ is 
surjective and 
\[
\sigma (Sf-Sg)\subset \sigma (f-g)
\]
holds for every pair $f$ and $g$ in $A$. We show that $S$ is linear 
and multiplicative. Let $\phi \in M_B$ be chosen arbitrarily. Then 
\[
\phi \circ S:A\to {\mathbb C},
\]
and 
\[
\phi \circ S(0)=0,
\]
and 
\[
\phi \circ S(f)-\phi \circ S(g)= 
\phi(Sf-Sg)\in \sigma (Sf-Sg)\subset \sigma (f-g)
\]
holds for every pair $f$ and $g$ in $A$. Thus by a theorem of 
Kowalski and S\l odkowski we have that $\phi \circ S$ is linear and 
multiplicative for every $\phi \in M_B$. Then conclusion follows immediately 
since $B$ is semi-simple.

We show that $T(0)=0$ if $a+b\ne 0$. Putting $f=g=0$ we have
\[
\sigma (aT(0)+bT(0))\subset \sigma (a\cdot 0 +b \cdot 0)=\{0\}.
\]
Thus we have $T(0)=0$ if $a+b\ne0$. 
\end{proof}

\section{A theorem of Moln\'ar and its generalizations}
On the other hand Moln\'ar \cite{m1} proved the following. 

\begin{theorem}(Moln\'ar)
Let ${\mathcal X}$ be a first countable compact Hausdorff space.
Suppose that $T$ is a map from $C({\mathcal X})$ onto itself such 
that the equality
\[
\sigma (TfTg)=\sigma (fg)
\]
holds for every pair $f$ and $g$ in $C({\mathcal X})$. Then there exist a 
continuous function $\eta:{\mathcal X}\to \{-1,1\}$ and a 
self-homeomorphism $\Phi$ 
on ${\mathcal X}$ such that the equality 
\[
Tf= \eta f\circ \Phi
\]
holds for every $f\in C({\mathcal X})$. In particular, $T$ is an algebra 
isomorphism if $T1=1$.
\end{theorem}
Motivated by the above theorems and others we may consider the 
following question: let $A$ and $B$ be Banach algebras and $p$ a polynomial 
of two variables. Suppose that $T$ is a map from $A$ into $B$ such that the 
inclusion
\[
\sigma (p(Tf,Tg))\subset \sigma (p(f,g))
\]
holds for every pair $f$ and $g$ in $A$. Does it follow that $T$ is 
linear and multiplicative? A theorem of Kowalski and S\l odkowski states that 
it is the case for $B={\mathbb C}$ and $p(z-w)=z-w$. On the other hand there 
several negative answers to the above too general question (see \cite{hmt1}). 
Even the polynomial $p$ need some restriction for a positive answer.

\begin{example}
Let $X$ be a compact Hausdorff space. For each $f\in C(X)$, put 
$\varepsilon_f =1$ or $-1$. Then the map $T$ from $C(X)$ into itself 
defined by 
\[
Tf=\varepsilon _f f,\quad f\in C(X)
\]
can be non-linear nor multiplicative but surjective. 
Put $p(z,w)=z^2+w^2$. Then the equality
\[
\sigma (p(Tf,Tg))=\sigma (p(f,g)), \quad f,g \in C(X)
\]
holds.
\end{example}

One of the reasonable questions may be as follows.

\noindent
{\bf Question}. {\it Let $A$ and $B$ be unital semi-simple commutative Banach 
algebras. Characterize the two-variable polynomials $p$ which satisfy 
the following property: a map $T$ from $A$ onto $B$ such that the 
equality
\[
\sigma (p(Tf,Tg))=\sigma (p(f,g)), \quad f,g \in A
\]
holds is an algebra isomorphism.}

A theorem of Moln\'ar gives a positive answer to the question, namely if
$A=B=C({\mathcal X})$, then $p(z,w)=zw$ is a desired polynomial. 
Theorem \ref{cks} states that for a Banach algebra $A$ and a semi-simple 
commutative Banach algebra $B$ $p(z,w)=az+bw$ is a desired polynomial.
If a type of a theorem of 
Kowalski and S\l odkowski for $p(z,w)=zw$ were true, positive results would 
follow for various Banach algebras with $p(z,w)=zw$. Unfortunately it is 
not the case; A modified {\it theorem} does not hold.
 On the other hand Moln\'ar \cite{m1} also proved a 
positive results for 
the Banach algebra of all bounded operators on an 
infinite-dimensional Hilbert space. 

Rao and Roy \cite{rr1} generalized a theorem of Moln\'ar for uniform algebras 
on the maximal ideal spaces 
and Hatori, Miura and Takagi \cite{hmt2} generalized for 
semi-simple commutative Banach algebras. For the case of uniform algebras, 
Hatori, Miura and Takagi \cite{hmt1} considered the equality of the range 
instead of that of the spectrum and show a generalization of a theorem 
of Moln\'ar. Luttman and Tonev \cite{lt} consider the equation for 
more smaller set; the peripheral range. Let $A$ be a uniform algebra on a 
compact Hausdorff space $X$. For $f\in A$, the peripheral range 
${\mathrm Ran}_{\pi}(f)$ for $f\in A$ is denoted by 
\[
{\mathrm Ran}_{\pi}(f)=\{z\in f(X): |z|=\|f\|_{\infty (X)}\}.
\]
Note that the peripheral range for uniform algebras coincides with the 
peripheral spectrum $\sigma_{\pi}(f)$;
\[
\sigma_{\pi}(f)=\{z\in \sigma (f):|z|=\rr (f)\}, 
\]
where $\rr (f)$ is the spectral radius. Luttman and Tonev proved the 
following.

\begin{theorem}(Luttman and Tonev) 
Let $A$ and $B$ be uniform algebras on compact Hausdorff spaces $X$ and $Y$ 
respectively. Suppose that $T$ is a map from $A$ onto $B$ such that the 
equality
\[
{\mathrm Ran}_{\pi}(TfTg)={\mathrm Ran}_{\pi}(fg)
\]
holds for every pair $f$ and $g$ in $A$. Then there exist a function $\eta :
M_B \to \{-1,1\}$ and a homeomorphism $\Phi$ from $M_B$ onto $M_A$ such 
that the equality
\[
\widehat{Tf}(y)=\eta (y)\hat f \circ \Phi (y), \quad y\in M_B
\]
holds for every $f\in A$, where $\hat \cdot$ denotes the Gelfand transform. 
In particular, $T$ is an algebra isomorphism if $T1=1$. 
\end{theorem}

\section{Main results}
\begin{theorem}\label{uniform}
Let $A$ and $B$ be uniform algebras on compact Hausdorff spaces $X$ and $Y$ 
respectively. 
Let $p(z,w)=zw +az+bw+ab$ be a polynomial. Suppose that $T$ is a map from 
$A$ onto $B$ such that the equality 
\[
{\mathrm Ran}_{\pi}(p(Tf,Tg))={\mathrm Ran}_{\pi}(p(f,g))
\]
holds for every pair $f$ and $g$ in $A$. 
Then we have the following.

(1) If $a\ne b$, then $T$ is an algebra isomorphism. Thus there exists an 
homeomorphism from $M_B$ onto $M_A$ such that 
\[
\widehat{Tf}(y)=\hat f\circ \Phi (y), \quad y\in M_B
\]
holds for every $f\in A$. 

(2)If $a=b$, then there exist a continuous map $\eta:M_B\to \{-1,1\}$ and a 
homeomorphism $\Phi$ from $M_B$ onto $M_A$ such that the equality 
\[
\widehat{Tf}(y)=\eta (y)\hat f\circ \Phi (y)+ a(\eta (y)-1), \quad y
\in M_B
\]
holds for every $f\in A$. 
\end{theorem}

The author does not know a similar result as Theorem \ref{uniform} 
holds for $p(z,w)=zw+az+bw+c$ ($ab\ne c$). 
In general for several polynomials a similar result as 
Theorem \ref{uniform} does not hold. For example let $p(z,w) = 
z^2+w^2$. Let $X$ be a disconnected compact Hausdorff space and 
$A=B=C(X)$. For each $f\in A$, $\eta_f$ is a map from $X$ into 
$\{-1,1\}$. Put a map $T$ from $A$ into $B$ by
\[
Tf=\eta_f f, \quad f\in A.
\]
Then we have 
\[
{\mathrm Ran}_{\pi}(p(Tf,Tg))={\mathrm Ran}_{\pi}(p(f,g))
\]
holds for every pair $f$ and $g$ in $A$. On the other hand $T$ may be 
surjective but non-linear nor multiplicative according to 
the choice of $\eta_f$.

\begin{proof}
Put a map $S:A \to B$ by 
\[
Sf=T(f-b)+b, \quad f\in A.
\]
By a simple calculation we see that $S(A)=B$ and 
\begin{equation}\label{1}
{\mathrm Ran}_{\pi}(S(f)(S(g)+c))=
{\mathrm Ran}_{\pi}(f(g+c))
\end{equation}
holds for every pair $f,g\in A$, where $c=a-b$. 

If $a=b$, then by a theorem of Luttman and Tonev \cite{lt} we see that 
there is a continuous function 
$\eta :M_B\to \{-1,1\}$ and a homeomorphism from $M_B$ onto $M_A$ such 
that 
\[
\widehat{Sf}(y)=\eta (y) \hat f\circ \Phi (y), \quad y\in M_B
\]
holds for every $f\in A$. It follows that 
\[
\widehat{Tf}(y)=\eta (y) \hat f\circ \Phi (y) + a(\eta (y) -1) 
\quad y\in M_B
\]
holds for every $f\in A$. 

Suppose that $a\ne b$. We show that $S$ is an isometric algebra 
isomorphism. First we show that $S$ is injective. To this end suppose that 
$Sf=Sg$. Then for every $h\in A$ we have
\begin{multline}
{\mathrm Ran}_{\pi}(fh)={\mathrm Ran}_{\pi}(S(f)(S(h-c)+c)) 
\\
= {\mathrm Ran}_{\pi}(S(g)(S(h-c)+c))
= {\mathrm Ran}_{\pi}(gh).
\end{multline}
Then by a routine argument applying peaking function argument we 
see that $f=g$. By putting $g=-c$ and $f\in A$ with $Sf=1$ in 
the equation \ref{1} we have
\[
\{0\}={\mathrm Ran}_{\pi}(f(-c+c))=
{\mathrm Ran}_{\pi}(S(-c)+c),
\]
so we have $S(-c)=-c$. Let $\lambda$ be an arbitrary complex number. 
Then we have
\begin{equation}\label{2}
\lambda {\mathrm Ran}_{\pi}(-cf)={\mathrm Ran}_{\pi}(\lambda (-c)f)=
{\mathrm Ran}_{\pi}(S(\lambda (-c))(S(f-c)+c)))
\end{equation}
and
\begin{multline}\label{3}
\lambda
{\mathrm Ran}_{\pi}(-cf)=\lambda {\mathrm Ran}_{\pi}
(S(-c)(S(f-c)+c))
\\
={\mathrm Ran}_{\pi}
(\lambda S(-c)(S(f-c)+c))=
{\mathrm Ran}_{\pi}
((-\lambda c)(S(f-c)+c))
\end{multline}
since $S(-c)=-c$. By a simple calculation
\[
B=\{S(f-c)+c: f\in A\}
\]
holds, and thus for every $G\in B$ we have
\[
{\mathrm Ran}_{\pi}(-\lambda c G)=
{\mathrm Ran}_{\pi}(S(-\lambda c)G)
\]
holds by the equations \ref{2} and \ref{3}. It follows that 
\[
-\lambda c = S(-\lambda c)
\]
holds and so
\[
\lambda =S(\lambda)
\]
holds for every complex number $\lambda$ since $c\ne 0$.

Next let $f\in A$. Then
\[
{\mathrm Ran}_{\pi}(f)={\mathrm Ran}_{\pi}
(S(1)(S(f-c)+c))={\mathrm Ran}_{\pi}(S(f-c)+c).
\]
We also see that
\[
{\mathrm Ran}_{\pi}(f)=
{\mathrm Ran}_{\pi}(S(f)(S(1-c)+c)=
{\mathrm Ran}_{\pi}(Sf)
\]
since $S(1-c)=1-c$.

Next let $P(A)$ be the set of all peaking functions in $A$. Then we see that 
\begin{equation}\label{4}
S(P(A))=P(B).
\end{equation}
Let $f\in P(A)$. Then $Tf\in P(B)$ since 
\[
\{1\}={\mathrm Ran}_{\pi}(f)={\mathrm Ran}_{\pi}(Sf).
\]
Note that $f$ is a peaking function if and only if 
${\mathrm Ran}_{\pi}(f)=\{1\}$. Thus we have that $S(P(A))\subset P(B)$ 
holds and the converse inclusion is proved in the same way since 
$S$ is a bijection. We also see by a simple calculation that 
\begin{equation}\label{5}
S(P(A)-c)+c=P(B).
\end{equation}

This does not prove  Theorem \ref{uniform}
 we can give the rest of the proof as in 
\cite{hmt1}, so we only sketch the rest of the proof. 

For $f\in P(A)$, put
\[
L_f=
\{x\in X:f(x)=1\}.
\]
Let $\Ch(A)$ be the set of all weak peak points for $A$. 
We denote for $x\in \Ch (A)$
\[
P_x(A)=\{f\in P(A):f(x)=1\}.
\]

{\bf Claim 1}. Let $f,g\in P(A)$. If $L_{Tf}\subset L_{Tg}$, 
then we have $L_f\subset L_g$. 

We show a proof. In the same way as in the proof of 
Lemma 2.2 in \cite{hmt1} we see that for every pair $f$ and $g$ in 
$P(A)$ the inclusion $L_f\subset L_g$ holds if and only if 
$1\in {\mathrm Ran}_{\pi}(ug)$ holds for every $u\in P(A)$ with 
$1\in {\mathrm Ran}_{\pi}(fu)$. Applying this and the equation 
\ref{5} we can prove Claim 1 in a way similar to the proof of 
Lemma 3.2 in \cite{hmt1}. 

{\bf Claim 2}. For every $y \in \Ch (B)$, there exists an $x\in 
\Ch (A)$ such that 
$S^{-1}(P_y(B))\subset P_x(A)$.

We show a proof. Let $f_1, \dots , f_n$ be a finite number of functions in 
$S^{-1}(P_y(B))$. We show that 
\[
\cap_{\j=1}^{n}L_{f_j}\ne \emptyset.
\]
Since $Sf_j\in P_y(B)$ we see that 
\[
\prod_{j=1}^nSf_j\in P_y(B).
\]
Since $S(A)=B$, there exists a $g\in A$ with $Sg=\prod_{j=1}^nSf_j$. 
Note that $g\in P(A)$ since $Sg\in P_y(B)$. We see that 
$L_{Sg}\subset L_{Sf_j}$ by the definition for every $j=1, \dots, n$. 
Then by Claim 1 we have that $L_g\subset L_{f_j}$ for every $j=1, 
\dots, n$, and so 
\[
L_g \subset\cap_{j=1}^nL_{f_j}.
\]
It follows that $\cap_{j=1}^nL_{f_j}\ne \emptyset$ since $g\in P(A)$ 
and so $L_g\ne \emptyset$. By the finite intersection property we see 
that 
\[
L=\cap_{f\in S^{-1}(P_y(B))}L_f \ne \emptyset.
\]
Since $L$ is a weak peak set for a uniform algebra $A$, there exists an 
$x\in L\cap\Ch (A)$. It follows that
\[
S^{-1}(P_y(B))\subset P_x(A).
\]
 
{\bf Claim 3}. For every $y\in \Ch (B)$, there exists a unique 
$x_y\in \Ch (A)$ such that 
\[
S(P_{x_y}(A))=P_y(B).
\]

We show a proof. Since $S^{-1}$ is a map from $B$ onto $A$ and the 
equality
\[
{\mathrm Ran}_{\pi}(S^{-1}(F)(S^{-1}(G)+c))=
{\mathrm Ran}_{\pi}(F(G+c)), \quad F,G\in B
\]
holds we can adapt a similar argument as in the proof of Claim 2 for 
$S^{-1}$ we see that for every $x\in \Ch (A)$ there exists a 
$y'\in \Ch (B)$ such that 
\[
S(P_x(A))\subset P_{y'}(B).
\]
Then by Claim 2 we see that for every $y\in \Ch (B)$ there exists an 
$x\in \Ch (A)$ and so $y'\in \Ch(B)$ such that 
\[
P_y(B)\subset S(P_x(A))\subset P_{y'}(B).
\]
It follows that $y=y'$ and the uniqueness of $x$ for $y\in \Ch(B)$.
We have proved Claim 3. 

We continue the proof of Theorem \ref{uniform}. Put a map $\phi :
\Ch (B) \to \Ch (A)$ by $\phi (y)=x_y$. Then in a similar way as 
in the proof of Theorem in \cite{hmt1} we see that the equality
\[
(S(f-c)+c)(y)=f\circ \phi (y), \quad y\in \Ch (B)
\]
holds for every $f\in A$. Substituting $f$ by $f-c$ we see that
\[
S(f)(y)=f\circ \phi (y), \quad y\in \Ch (B).
\]
It follows that $S$ is an algebra isomorphism from $A$ onto $B$. 
Thus by the routine argument of commutative Banach algebras 
we see that there exist a homeomorphism $\Phi$ from $M_B$ 
onto $M_A$ such that the equality
\[
\widehat{S(f)}(y)=\hat f\circ \Phi (y),\quad y\in M_B
\]
holds for every $f\in A$. Then by the definition of $S$ we see by a 
simple calculation that the equality
\[
\widehat{Tf}(y)=\hat f\circ \Phi (y), \quad y\in M_B
\]
holds for every $f\in A$.
\end{proof}

\begin{theorem}\label{semi}
Let $A$ be a unital semi-simple commutative 
Banach algebra and $B$ a unital commutative Banach algebra. Put 
$p(z,w)=zw+az+bw+c$, where $a,b$ and $c$ are coefficients. Suppose that
$T$ is a map from $A$ onto $B$ such that the equality
\[
\sigma (p(Tf,Tg))=\sigma (p(f,g))
\]
holds for every pair $f$ and $g$ in $A$. Then we have the following.

(1) If $a\ne b$, then $T$ is an algebra isomorphism. Thus there 
exists a homeomorphism from $M_B$ onto $M_A$ such that the equality
\[
\widehat{Tf}(y)=\hat f\circ \Phi (y), \quad y\in M_B
\]
holds for every $f\in A$.

(2) If $a=b$, then there exist a map $\eta:
M_B\to \{-1,1\}$ and a homeomorphism $\Phi$ from $M_B$ onto $M_A$ such that 
the equality
\[
\widehat{Tf}(y)=\eta (y)\hat f\circ \Phi (y) + a(\eta (y)-1), \quad y\in M_B
\]
holds for every $f\in A$. 

\noindent
In any case we have that $B$ is semi-simple and $A$ is algebraically 
isomorphic to $B$.
\end{theorem}

\begin{proof}
We consider the case where $B$ is semi-simple. 
(The general case follows from the case where $B$ is semi-simple. 
Consider the Gelfand transform $\Gamma$ of $B$. Then the composition map
$\Gamma \circ T$ is a map from $A$ onto the Gelfand transform $\hat B$ of 
$B$. Then by the first part we see that $\Gamma \circ T$ is 
injective, which will follow that $\Gamma$ is injective. 
Thus we see that $B$ is semi-simple and we can deduce the case where 
$B$ is semi-simple.) Put a map $S:A\to B$ by
\[
S(f)=T(f-b)+b, \quad f\in A.
\]
Then by a simple calculation we see that $S(A)=B$ and the equality
\[
\sigma (f(g+c))=\sigma(S(f)(S(g)+c)), \quad f,g \in A
\]
holds, where $c=a-b$.

If $a=b$, then by a proof of Theorem 3.2 in \cite{hmt2} there exist a 
continuous function $\eta:M_B\to \{-1,1\}$ and a homeomorphism 
$\Phi$ from $M_B$ onto $M_A$ such that the equality
\[
\widehat{Sf}(y)=\eta (y)\hat f\circ \Phi (y), \quad y\in M_B
\]
holds for every $f\in A$. It follows that 
\[
\widehat{Tf}(y)=\eta (y) \hat f\circ \Phi (y) + a(\eta (y) -1) 
\quad y\in M_B
\]
holds for every $f\in A$. 

Suppose that $a\ne b$. Then by the same way as in the proof of Theorem 
\ref{uniform} we see that $S(-c)=-c$ and the equality
\[
S\lambda =\lambda
\]
holds for every complex number $\lambda$.

{\bf Claim 1}. For every $f\in A^{-1}$, the equality $S(f)(S(f^{-1}-c)+c)=1$. 

We show a proof. Since
\[
\{1\}=\sigma(ff^{-1})=\sigma (S(f)(S(f^{-1}-c)+c)
\]
we have
\[
S(f)(S(f^{-1}-c)+c)=1
\]
since $B$ is semi-simple.
We denote the uniform closure of $\hat A$ in $C(M_A)$ by $\mcl (A)$, 
where $C(M_A)$ is the algebra of all complex-valued continuous functions 
on $M_A$. Note that the maximal ideal space of $\mcl (A)$ coincides with 
$M_A$. In the following the Gelfand transformation of $f$ in $A$ and 
$\mcl (A)$ is denoted also by $f$ for simplicity.

{\bf Claim 2}. Let $\{f_m\}$ be a sequence in $A^{-1}$ and 
$f\in C(M_A)$ such that 
\[
\|f_m-f\|_{\infty (M_A)} \to 0
\]
as $m\to \infty$. Then $\{Sf_m\}$ is a Cauchy sequence in $B$ with respect to 
the uniform norm on $M_B$ and the uniform limit $\lim Sf_m$ is an 
invertible function in $\mcl(B)$. 

We show a proof of Claim 2. We may assume that there exists a positive 
integer $K$  with the inequality
\[
\frac{1}{K}<|f_m(x)|<K, \quad x\in M_A
\]
holds for every positive integer $m$. Note that 
\[
\frac{1}{K}<|Sf_m(y)|<K, \quad y\in M_B
\]
holds for every positive integer $m$ since 
\[
\sigma(f_m)=\sigma(Sf_m(S(1-c)+c))=\sigma(Sf_m)
\]
holds. 
Then by a simple calculation we see 
that for every positive $\varepsilon$, there exists a positive integer $N$ 
such that the inequality
\[
\left|\frac{f_n(x)}{f_m(x)}-1\right|<\varepsilon, \quad x\in M_A
\]
holds for every $m,n>N$. Since $Sf_m(S(f_m^{-1}-c)+c)=1$ we see that 
\[
\sigma (f_nf_m^{-1})=\sigma (Sf_n(S(f_m^{-1}-c)+c)=
\sigma (Sf_n(Sf_m)^{-1}),
\]
so the inequality
\[
\left|\frac{Sf_n(y)}{Sf_m(y)}-1\right|<\varepsilon, \quad y\in M_B
\]
holds for every $m,n >N$. Thus we see that
\[
\|Sf_n-Sf_m\|_{\infty (M_B)}\le \|Sf_n\|_{\infty (M_B)}
\|\frac{Sf_n}{Sf_m}-1\|_{\infty (M_B)} \le K\varepsilon
\]
holds for every $m,n > N$, so $\{Sf_m\}$ is a Cauchy sequence with 
respect to the uniform norm and 
\[
\frac{1}{K}\le |\lim Sf_m|\le K
\]
on $M_B$, so $\lim Sf_m$ is invertible in $\mcl (B)$ since the maximal 
ideal space of $\mcl (B)$ coincides with $M_B$.
We have proved Claim 2. 

{\bf Claim 3}. Then map $S$ is extended to an injective map $\bar S$ 
from $A\cup (\mcl(A))^{-1}$ onto $B\cup (\mcl (B))^{-1}$ such that 
the equality
\[
{\mathrm Ran}(\bar Sf(\bar Sg+c))={\mathrm Ran}(f(g+c))
\]
holds for every pair $f$ and $g$ in $A\cup (\mcl (A))^{-1}$.

We show a proof. Let $f\in (\mcl (A))^{-1}$. Note that 
\[
(\mcl (A))^{-1}=\{f\in \mcl (A):0\not\in f(M_A)\}
\]
since the maximal ideal space of $\mcl (A)$ coincides with $M_A$. 
Then there exists a sequence $\{f_m\}$ in $A$ with 
\[
\|f_m-f\|_{\infty (M_A)}\to 0
\]
as $m\to \infty$. We may assume that $f_m\in A^{-1}$. Then by 
Claim 2 we see that the uniform limit $\lim Sf_m$ exists and it is 
easy to see that the limit does not depend on the choice of 
a sequence $\{f_m\}$ which converges to $f$. Put $\bar Sf=
\lim Sf_m$. Then by Claim 2 we see that $\bar S f \in (\mcl (B))^{-1}$. 
In this way we can define $\bar S$ from $A\cup (\mcl (A))^{-1}$ 
into $B\cup (\mcl (B))^{-1}$. By some calculation 
we see that 
\[
{\mathrm Ran}(\bar Sf(\bar Sg +c))=
{\mathrm Ran}(f(g+c)), \quad f,g \in A\cup (\mcl (A))^{-1}
\]
holds. We also see in the same way as in the 
proof of Claims 3 and 4 in \cite{hmt1} that $\bar S$ is a bijection. 

This does not prove the theorem, but the rest of the proof is 
similar to that of a proof of Theorem 3.2 applying a similar way as in 
the proof of Theorem \ref{uniform}. We omit a precise proof. 
\end{proof}

\end{document}